\documentclass[12pt]{article}
\usepackage{latexsym}
\usepackage{amsfonts}
\usepackage{times}

\newcommand{\pnt}{$\!\!\! \bf . \;$}
\newtheorem{lemma}{Lemma}
\newtheorem{corollary}{Corollary}[lemma]
\newtheorem{theorem}[lemma]{Theorem}

\title{The Circuit Polynomial of the Restricted Rooted Product
$G({\it\Gamma})$ of Graphs with a Bipartite Core $G$} 
\author{Vladimir R. Rosenfeld\\
{\small \sl Institute of Evolution, University of Haifa, Mount Carmel,
Haifa 31905, Israel}\\
{\small \sl E-mail: vladimir@research.haifa.ac.il}} 

\date{ }
\begin{document}
\setcounter{page}{1}

\maketitle

\vspace*{2 cm}

\begin{abstract}
As an instance of the $B$-polynomial, the circuit, or cycle, polynomial
$P(G({\it\Gamma}); {\rm w})$ of the generalized rooted product
$G({\it\Gamma})$ of graphs was studied by Farrell and Rosenfeld [19]
and Rosenfeld and Diudea [20]. In both cases, the rooted product
$G({\it\Gamma})$ was considered without any restrictions on graphs $G$ and
${\it\Gamma}$. Herein, we present a new general result and its corollaries
concerning the case when the core graph $G$ is restricted to be bipartite.
The last instance of $G({\it\Gamma})$, as well as all its predecessors
[19, 20], can find chemical applications.
 
\end{abstract}

\section{Introduction}

\hspace*{6mm}
As it was phrased in Scientific American, today is a nanotechnology gold
rush. To the mathematical reader, it may be pleasant to recognize that
mathematics (and, in particular, the theory of graph spectra) is also
seeking to make its contribution to this interdisciplinary area
[1--20]. The present paper continues the previous work [19, 20] along
these lines.

From the graph-theoretical point of view, there exists a number of molecular 
structures (or graphs) of high-tech interest that can be generated using
the graphical construction called, in the mathematical literature, the
generalized rooted product of graphs (see [19, 20, 29, 30]). Here,
first of all, of special practical significance are dendrimers
[1--8], for which some authors also use other names (e.g., "bundled
structures", as in [9]).

As well as any other type of products relevant to chemical objects, the
rooted product of molecular graphs contains combinatorial information that
is useful for the (theoretical) chemist. Algebraically, this information
can be represented in the form of respective graph polynomials of the
molecular graph in question. One type of such polynomials (namely, the
circuit polynomial) will be considered by us throughout this paper.

An algorithmic role of the results that will be discussed below can be
clarified through recalling, as an example, the following situation.
Creating complex dendritic molecules also foresees that the researcher
would have suitable quantum-chemical (or any other) methods for estimating
the energy levels of respective large molecules. A common tack is to
utilize for this the simpler solutions that had been obtained for
dendrimer's core and monodendrons. Namely, this sort of calculation
propose mathematical results that will follow below. However, we want to
specially \underline{stress that} that we shall, in fact, deal with the
derivation of universal mathematical relationships for the circuit
polynomial which are not based on any approximation might follow from
quantum-chemical or other calculational methods that can adapt them. 

Now that we have briefly described the applied background of the paper,
we need to focus upon a rigorous mathematical exposition of our specific
task.

\section{Preliminaries}

\hspace*{6mm}
We should start with giving some notions from the theory of graph 
polynomials.

\subsection{The {\bf ${F}$}- and {\bf ${B}$}-polynomials of a graph} 

\hspace*{6mm}
The graphs considered here are finite, may be directed, weighted; and may
contain self-loops, i.e., finite directed or undirected weighted
pseudographs. A general class of graph polynomials was introduced in
Farrell [21]. These are called $F$-polynomials and are defined as
follows. Let $G$ be a graph and $F$ a family of connected subgraphs of
$G$. An $F$-{\em cover} of $G$ is a spanning subgraph of
$G$, in which every component is a member of $F$. Let us associate with
each member $\alpha$ of $F$ an intermediate or weight $w_{\alpha}$. The
weight of a cover $C$ denoted by $w(C)$, is the product of the weights of
its components. Then the $F$-{\em polynomial} is
\begin{equation}
\label{1}
F(G;{\bf w})=\sum w(C)\,,
\end{equation}
where the summation is taken over all the $F$-covers of $G$, and where
${\bf w}$ is a vector of the indeterminates $w_{\alpha}$.

Throughout this paper, we denote the vertex (or node) set of $G$ by $V(G)$
and assume that $|V(G)|=p$, unless otherwise specified. Also, if G is
labeled, we associate with the $i$-th vertex of $G$ the special weight
$x_{i}+b_{i}$\, $(1\leq i\leq p)$, where $x_{i}$ is an
indeterminate and $b_{i}$ is the the sum of weights of all loops, if any,
lying in a vertex $i$ (see [22, 23]). We use the notation $F(G;{\bf
x})$, for $F(G;{\bf w})$, when all the variables, except the $x_{i}$'s, are
replaced by $1$'s. If we replace all $x_{i}$'s, in $F(G;{\bf x})$, with
the single variable $x$, then the resulting polynomial in $x$ will be
denoted by $F(G;x)$, and called the {\em simple} polynomial of $G$.

If every nonnode member of $F$ consists of exactly one block, then we
call the corresponding class of $F$-polynomials, {\em block
polynomials}; or $B$-{\em polynomials}, for short. We
then write $B(G;{\bf w})$ for $F(G;{\bf w})$, in order to indicate this
property of the members of $F$. Notice that if we take $F$ to be a family
of cycles, then every nonnode member of $F$ is a block. This is also true
when $F$ is the family of cliques. Therefore both the circuit (or
cycle) polynomial and clique polynomial (see [19]) are examples of
block polynomials. We therefore classify all the special circuit
polynomials, for example the matching, characteristic and permanental
polynomials (see [10--25; 28--30]), as $B$-polynomials.  

It should be observed that the families which give rise to $B$-polynomials
consist of graphs which are characterized by the number of vertices.
Therefore, when general weights are to be assigned to the members of $F$,
it is sufficient to associate with each member, of $F$, with $n$ vertices
the weight $w_{n}$. The resulting $B$-polynomial would therefore contain
monomials which totally describe the covers. In this general  
$F$-polynomial, the vector of weights is ${\bf w}=(w_{1},w_{2},\ldots, 
w_{p})$. Observe that if $F$ is the family of stars or paths, then
every member of $F$ is characterized by the number of nodes. However,
stars and paths are not blocks and so do not give rise to $B$-polynomials.

The stimulus to investigate the $B$-polynomials stems from the fact that
they are often encountered in many problems in Mathematics, as well as in
various applications outside of Mathematics. It is interesting to know
about mutual and hereditary relations among different graph polynomials.
For instance, the matching polynomial is a generalization of the so-called
acyclic polynomial, which was defined independently (see [19]). The same
matching polynomial yields, under certain substitutions, the chromatic
polynomial for certain classes of graphs, and also a whole group of its
relatives (see [19]) as well. The classical rook polynomial (see [19])
is yet another relative of the matching polynomial. 
                                                   
Notice that the most general $F$-polynomial is the subgraph polynomial
(see [19]), since it enables us to derive, in principle, any other
$F$-polynomial. However, the subgraph polynomial is not a $B$-polynomial.
So, there exist other classes of $F$-polynomials; e.g., see [19],
wherein the so-called {\em articulation node polynomials} (or
$A$-{\em polynomials}, for short) are introduced.

Now we shall specially consider some instances of the circuit polynomial.

\subsection{The circuit (cycle) polynomial of a graph}

\hspace*{6mm}
The circuit (cycle) polynomial $C(U;{\bf w})$ of an undirected graph $U$
was introduced by Farrell [24] (see [25]). The notion of this polynomial
was generalized in [19, 20] for an arbitrary graph $G$. Herein, we shall
give the third definition of it, which is, however, tantamount to that of
[19, 20], where the circuit polynomial was considered in quite a different
way, as a specific case of the $F$-polynomial. In order to indicate the
distinction between the original Farrell's polynomial $C(U;{\bf w})$ and
the one that will be used in the present paper, we shall denote the latter
by $P(U;{\bf x};{\bf w})$, where ${\bf x}=(x_1,x_2,\ldots,x_p)$ is, by
analogy with ${\bf w}$, the vector of indeterminates (see [20]).  

Many properties of the polynomial in question can be considered from the
matrix-theoretical standpoint. Let $A=\{a_{ij}\}_{i,j=1}^p$ be the
adjacency matrix of a graph $G$. Let further $A^{*}=\{a_{ij}^{*}\}=(A+X)$
be an auxiliary matrix, where $X$ is a diagonal matrix, whose on-diagonal
entries are indeterminates $x_{1},x_{2},\ldots,x_{p}$, consecutively.
One can define {\em circuit} ({\em cycle}) {\em polynomial} $P(G;{\rm
x};{\rm w})$ {\em of a graph} $G$ as follows

\begin{equation}
\label{2}
P(G;{\bf x};{\bf w}):=\sum_{\sigma\in S_{p}}\prod_{i=1}^{p}a_{i\sigma 
i}^{*}w_{i}^{\omega_{i}(\sigma)},
\end{equation}
where $a_{i\sigma i}^{*}$ is the respective entry of $A^{*}$;
$\omega_{i}(\sigma)$ is the number of cycles of length $i$ in a
permutation $\sigma$; and the sum ranges over all the $p!$ permutations 
$\sigma$ of a symmetric group $S_{p}$. (Recall that $\sigma i=j$ is the
image of an index $i$, obtained under the action of a permutation $\sigma$
on the set $I=\{1,2,\ldots,p\}$ of vertex indices; $i,j\in I$.)

The polynomial $C(U;{\bf w})$ of an undirected graph $U$, introduced by
Farrell (see [24, 25]), is a specific case of $P(U;{\bf x};{\bf w})$,
viz.:
\begin{equation}
\label{3}
C(U;{\bf w})=P(U;{\bf x};{\bf w})\left |_{x_{i}=1;\, b_{i}=0;\,
w_{j}\rightarrow\frac{w_{j}}{2}}\right.\hspace{6mm}(1\leq i\leq p;3\leq
j\leq p),
\end{equation}
where $w_{j}\rightarrow\frac{w_{j}}{2}$ denotes the substitution of
$\frac{w_{j}}{2}$ for $w_{j}$.

In general, if a graph $G$ is directed (2) does not hold. We can mention
in passing one ready result for $C(U;{\bf w})$, connected to the
enumerative theory of P\'{o}lya (see [26, 27]). Let $K_{p}$ be a
complete
graph with $p$ vertices. Then
\begin{equation}
\label{4}
C(K_{p};{\bf w})=p!Z(S_{p};V;w_{1},w_{2}\ldots,w_{p}),
\end{equation}
where $Z(S_{p};V;{\bf w})$ is the cycle indicator of a symmetric group
$S_{p}$ faithfully acting on a vertex set $V=V(K_p)$ $(|V|=p)$ of
$K_p$ (see [26, 27]).  

The circuit polynomial $P(G;{\bf x};{\bf w})$ has, as its specific cases,
the {\em generalized permanental polynomial} $\phi^{+}(G;{\rm x})$, 
{\em generalized characteristic polynomial} $\phi^{-}(G;{\rm x})$, and
two {\em generalized matching polynomials} 
$\alpha^{+}(G;{\rm x})$ and $\alpha^{-}(G;{\rm x})$ [22, 23]; viz.:
\begin{equation}
\label{5}
\phi^{+}(G;{\bf x})={\rm per}(A+X)=P(G;{\bf x};{\bf 
w})\left |_{w_{i}=1}\right.\hspace{6mm}(1\leq i\leq p),
\end{equation}
\begin{equation}
\label{6}
\phi^{-}(G;{\bf x})={\rm det}(A-X)=P(G;{\bf x};{\bf 
w})\left |_{w_{i}=-1}\right.\hspace{6mm}(1\leq i\leq p),
\end{equation}
\begin{equation}
\label{7}
\alpha^{+}(G;{\bf x})=P(G;{\bf x};{\bf
w})\left |_{w_{1}=w_{2}=1;\,w_{i}=0}\right.\hspace{6mm}(3\leq i\leq p),
\end{equation}
\begin{equation}
\label{8}
\alpha^{-}(G;{\bf x})=P(G;{\bf x};{\bf
w})\left |_{w_{1}=w_{2}=-1;\,w_{i}=0}\right.\hspace{6mm}(3\leq i\leq p).
\end{equation}

It is worth noting that the weight $b_{i}$ $(1\leq i\leq p)$ (see [19,
20]) of a self-loop lying in a vertex $i$ of $G$ is thereby equal to
an entry $a_{ii}$, of $A$, in case of $\phi^{+}(G;{\bf x})$ and
$\alpha^{+}(G;{\bf x})$; however, $b_{i}$ is equal to $-a_{ii}$ in case of
$\phi^{-}(G;{\bf x})$ and $\alpha^{-}(G;{\bf x})$.

Apparently, there are more possibilities to devise other such
polynomials with the adjective "generalized"; herein, we shall confine
ourselves only with the above instances, to avoid any confusion. Recall
the {\em simple circuit polynomial} $P(G;x)$ is a one-variable case of
it, with an italicized $x$ in a lieu of \textbf{x}, viz.: 

\begin{equation} 
\label{9}
P(G;x)=P(G;{\bf x})\left |_{x_{i}=x}\right.\hspace{6mm}(1\leq i\leq p),
\end{equation}
while the variables $w_{1},w_{2},\ldots,w_{p}$ may or may not be reduced
(it depends on the context).

The notation $P(G;{\bf x};{\bf w})$ or any its reduced-variable form will
hereafter stand for every possible instance of it at once; the
reader can reinterpret any of general solutions for any specific circuit
polynomial that he/she needs in---the permanental, characteristic,
matching. (Moreover, some other $B$-polynomials can have the same
properties, e.g., the clique polynomial; see [19].) However, one
would recognize that the most studied and widely used instance of the
circuit polynomial is the {\em simple characteristic polynomial}
$\phi^{-}(G;x)$ (see [10--18; 23, 24; 28--31], where [28, 29]
and [11] (for chemists) are the main world's monographs on the
subject). 

In order to proceed, we need to consider now some kinds of operations on
graphs.

\subsection{Some products of graphs}

\hspace*{6mm}
Let $(G,u)$ and $(H,v)$ be two graphs rooted at node $u$ and $v$,
respectively. We {\em attach} $G$ to $H$ (or $H$ to $G$) by identifying
node $u$ of $G$ with node $v$ of $H$. Nodes $u$ and $v$ are called
{\em nodes of attachment}. The node formed by identification is called
the {\em coalescence node}. The resulting graph ${G\circ H}$ is called 
the {\em coalescence of} ${G}$ {\em and} ${H}$. 

Now consider a family 
$\{(U_{1},u_{1}),(U_{2},u_{2}),\ldots,(U_{t},u_{t})\}$ of not necessary
distinct graphs with roots $u_{1},u_{2},\ldots,u_{t}$, respectively. We
term a connected graph $U_{1}\circ U_{2}\circ\cdots U_{t}$ the
{\em multiple coalescence of} ${U_{1},U_{2},\ldots,U_{t}}$ 
\underline{provided that nodes {$\mathit {u_{1},u_{2},\ldots,u_{t}}$} 
are identified to reform the} \underline{coalescence node {$\mathit
{r}$}}. We shall use ${U^{|q|}}$ to denote a $q$-fold coalescence of $q$ 
isomorphic copies of a graph $U$; in the same way, we shall use
${G\circ H^{|s|}}$ to denote the multiple coalescence of $G$ and $s$
copies of $H$, wherein all coalesced graphs have just one cutnode $r$
in common; etc.. 

The above operation {\boldmath ${\circ}$} is associative; in other words,
it can be met as a generating operation in some semigroups of graphs. As a
case in point, pick the set ${\mathcal {U}}=\{(U_{j},u_{j})\}_{j=1}^ 
{\infty}$ of all unicomponental graphs; obviously, a pair $({\mathcal
{U}};\circ)$ is an infinite commutative monoid of graphs, wherein the
unity is represented by a one-vertex graph $K_{1}$. 

Let $G$ be a graph with $p$ nodes and $\Gamma =\{
H_{1},H_{2},\ldots,H_{n}\}$ a family of rooted graphs. Then the graph
formed by attaching $H_{k}$ to the $k$-th $(1\leq k\leq p)$ node of $G$ is
called the {\em generalized rooted product} (see [19, 20]) and
is denoted by ${G({\it\Gamma})}$; $G$ itself is called the
{\em core} of $G({\it\Gamma})$. If each member of ${\it\Gamma}$ is
isomorphic to the rooted graph $H$, then the graph $G({\it\Gamma})$ is
denoted by ${G(H)}$ [19, 20, 29, 30]. Furthermore, if $H$ is
an edge (a twig), then the resulting graph is called a {\em thistle} or
{\em equible graph} (see [19]). 

Herein, we should make it our first business to treat some specific cases of 
the rooted product. Let a core $G$ be a bipartite graph $T$, whose parts
have $p_{1}$ and $p_{2}$ vertices, respectively; $p_{1}\geq p_{2}$ and 
$p_{1}+p_{2}=p$. One can attach to every vertex in the first part of $T$
an isomorphic copy of a graph $H_{1}$ and to every other vertex, in
$B$, an isomorph of another graph $H_{2}$ $(H_{1},H_{2}\in
{\it\Gamma})$; we shall locally call the resulting graph 
$T({\it\Gamma})$ the {\em restricted rooted product of graphs}
$T$ {\em and} ${\it\Gamma}$. It is worth specially noting an instance
of it in which one type of graphs (either $H_{1}$ or $H_{2}$) is simply
a one-vertex graph $K_{1}$ (with or without self-loops); that is, thus,
in either part of $T$ no attachment is made. Since the last case
promises a nice operation for constructing more complex graphs of
practical interest, we shall supply some relevant information about
them.

An interesting generalization of the rooted product are the
${F}$-{\em graphs} [31], which are consecutively iterated rooted
products defined as follows: $F^{0}=K_{1},F^{1}=G,F^{2}=G(H),
\ldots,F^{s+1}=F^{s}(H)$\, $(s\geq 1)$. Another interesting example of 
rooted product is the family of {\em dendrimers} ${D^{k}}$\,
$(k\geq 0)$\, (see [7, 8, 19, 20]), defined as follows: $D^{0}=K_{1},
D^{1}=G,D^{2}$ is the rooted product of $G$ and $H$, in which
\underline{some} attachments of $H$ are not made, i.e., $H$ need not be
attached to all nodes of $G$. In general, $D^{s+1}$\, $(s\geq 1)$ is
constructed from $D^{s}$; and the number of copies attached to $D^{s}$
obeys some fixed generation law. The dendrimers, in particular, imitate
molecular structures, bearing the same name [1--8; 19, 20]. They are of
practical significance [1--6]. This has lent impetus to the 
investigations in this paper as well.

A {\em monodendron} $M$ is a maximal connected subgraph of
a dendrimer $D$ that shares only the coalescence node $r$ with a core $G$;
in other words, it is a maximal (hyper)branch of $D$. Being a dendrimer in
its own right, $M$ has, however, two peculiarities. First, its core $G$ is 
played by the same (weighted) graph $H$ that is a structural repeating
unit of branches $(G=H)$. Second, a core $G$\,(or $H$) of $M$ possesses
the root (node $r$), which is not a feature of all dendrimers. Owing to
its root $r$, the entire monodendron $M$ can be made to serve the function 
of a new structural unit (instead of $H$) for constructing the higher
dendrimers. Moreover, as well as any other dendrimer, $M$ can serve as a
hypercore in the same procedure, in lieu of a simple core $G=D^{1}$. As an
instance of ${\mathcal D}=\{ D^{j}\}_{j=0}^{\infty}$, the monodendron
series ${\mathcal M}=\{ M^{j}\}_{j=0}^{\infty}$ is defined as 
follows: $M^{0}=K_{1},M^{1}=G=H$ and $M^{k}$\, $(k\geq 2)$ is constructed
by analogy with $D^{k}$ above.  

Let (a copy of) $H$ invariably make $d+1$\, $(d<p(H))$ attachments inside
a dendrimer $D$. Of this amount, 1 attachment is to hold the root of $H$
itself while the other $d$ are to hold the roots of all its incident
neighbors in $D$. The number $d$ is called a {\em progressive degree
of} $H$ ({\em cf} [8]). A dendrimer is said to be
{\em homogeneous} if all its monodendrons are equivalent and all 
prescribed attachments within it are made ({\em cf} [8]). By definition,
all dendrimers that we consider herein are homogeneous. 

A monodendron $M^{j}$\, $(j\geq 1)$ contains $1+d+d^{2}+\cdots
+d^{j-1}=(d^{j}-1)/(d-1)$ isomorphic copies of $H$ therein; they are
lying in concentric {\em layers} ({\em tiers}). This distribution
correlates with their distance from a core $G$; all copies of $H$ that are
built into one and the same tier are spaced at the same distance from $G$.

It is convenient to begin numbering the layers in a monodendron $M$ from
its core $G=H$ (thus receiving the ordinal 1). So, the number of layers in
$M^{j}$\, $(j\geq 0)$ equals to $j$ itself and the $k$-th layer contains
$d^{k-1}$ isomorphic copies of $H$ (and 0 under $k=0$); the number of
nodes, in the $j$-th (uttermost) layer, that can be used for further
attachments is $d^{j}$.

We come now to an important remark. The matter is that the above procedure
that successively generates all monodendrons $M^{j}$ of a series
{$\mathcal M$} is unambiguous; it always reproduces one and the same
monodendron $M^{j}$ with a given number $j$ of tiers. So, it is impossible
to produce instead of $M^{j}$ any other homogeneous monodendron with $j$
layers. Hence it follows that the number $j\geq 0$ of layers
\underline{uniquely characterizes} a homogeneous monodendron
$M^{j}$ in a specific series {$\mathcal M$}.                                                               

In a general way, a monodendron $M^{j}$ and $d^{j}$ isomorphic copies of a
monodendron $M^{k}$\, $(j,k\geq 0)$, when used as $G$ and $H$,
respectively, afford a monodendron $M^{j+k}$ (see above). We shall use the
notation $M^{j}\star M^{k}=M^{j+k}$ to denote this. The binary operation
{\boldmath $\star$} is obviously commutative and associative: $M^{j}\star
M^{k}=M^{k}\star M^{j}$ and $M^{j}\star (M^{k}\star M^{l})=(M^{j}\star
M^{k})\star M^{l}$ $(j,k,l\geq 0)$, which can readily be verified,
recalling that the number of tiers in the resulting monodendron uniquely
characterizes it. Since $M^{0}=K_{1}$ acts as the identity, we can at once
conclude that $({\mathcal M};\star)$ is an infinite commutative monoid
isomorphic to the additive monoid ({\boldmath $N$};+) of all nonnegative
integers. One can simply say that {$\mathcal M$} is a monoid (without
indicating its operation) and also adopt the multiplicative notation
$M^{j}M^{k}$ for $M^{j}\star M^{k}$ and in any similar case. The said of
{$\mathcal M$} herein resembles two earlier-studied situations [31, 32]
(see [33]) in every essential detail. 

Now we need to consider some known results that will be used below.

\subsection{Basic results}

\hspace*{6mm}
We begin with a previous result (see Lemma 5 in [19]), rewritten
here as follows:
\begin{lemma}\pnt
Let $G\circ H$ be the graph formed by attaching a graph $G$ to a graph
$H$, and let $r$ be the resulting coalescence node. Then
\begin{equation}
\label{10}
B(G\circ H;{\bf x};{\bf w})=B(H_{-r};{\bf x};{\bf w})\left [B(G;{\bf
x};{\bf w})\left |_{x_r\rightarrow\frac{B(H^{\triangle};{\bf x};{\bf
w})}{B(H_{-r};{\bf x};{\bf w})}}\right.\right ]\, ,
\end{equation}
where $H_{-r}$ is a graph $H$ less its root $r$; and $H^{\triangle}$ is
the graph $H$ less its self-loops lying in the vertex $r$.
\end{lemma}

An important result is the following statement (see Theorem 2 in [19]):
\begin{theorem}\pnt
Let $G$ be a graph with $p$ vertices and 
${\it\Gamma}=\{H_1,H_2,\ldots,H_p\}$ be a family of rooted graphs. Then
\begin{equation}
\label{11}
B(G({\it\Gamma});{\bf x};{\bf w})=\left [\prod_{i=1}^{p}B(L_{i};{\bf 
x};{\bf w}\right ]\left [B(G;{\bf x};{\bf w})\left
|_{x_i\rightarrow\frac{B(H^{\triangle};{\bf x};{\bf
w})}{B(L_{i};{\bf x};{\bf w})}}\right.\right ]\hspace{6mm}(1\leq i\leq
p),
\end{equation}
where $L_{i}$ $(1\leq i\leq p)$ is the graph $H_{i}$ with its root removed
(i.e., $H_{i}-r_{i}$); and $H_{i}^{\triangle}$ is the graph $H_{i}$ with
all the self-loops at the root removed.
\end{theorem}

In the case of the simple rooted product, one can derive the following
corollary of Theorem~2 for the simple $B$-polynomial (see Corollary 2.1 in
[19]), viz.:
\begin{corollary}\pnt
Let $G$ and $H$ be rooted graphs. Let $G(H)$ be the graph obtained by
attaching an isomorph of $H$ to each of the $p$ nodes of $G$. Then
\begin{equation}
\label{12} 
B(G(H);x)=[B(H_{-r};x)]^{p}\left [B(G;x)\left
|_{x\rightarrow\frac{B(H^{\triangle};x)}{B(H_{-r};x)}}\right.\right ],
\end{equation}
where $H^{\triangle}$ is the graph $H$ with all loops at its root $r$
removed.
\end{corollary}
Here, we should note earlier specific versions of Lemma 2 for the
characteristic [30, 29] and matching [29] polynomials (wherein only
unweighted graphs have been treated).

Now recall that any simple $B$-polynomial, such as $B(G;x)$ in (11), can
be expanded in powers of $x$; therefore we can write down it as
\begin{equation}
\label{13}
B(G;x)=\gamma_{0}x^{p}+\gamma_{1}x^{p-1}+\ldots
+\gamma_{p}x^{0}\hspace{6mm}(\gamma_0=1)\, .
\end{equation}

Owing to (13), we can give herein a new version of Lemma 3 (see Lemma 3 in
[20]); viz.:
\begin{lemma}\pnt
Let $G$ and $H$ be rooted graphs. Let $G(H)$ be the rooted product of $G$
and $H$, as above. Then
\begin{equation}
\label{14}
B(G(H);x)=\sum_{g=0}^{p}\gamma_{g}[B(H^{\triangle};x)]^{p-g}[B(H_{-r};x)]^{g} 
\, ,
\end{equation}
where $H^{\triangle}$ is the same as above.
\end{lemma}

The next quotation (Corollary 2.2 from [19]) appears herein as follows:
\begin{lemma}\pnt
Let $G$ and $H$ be rooted graphs. Let $G(H)$ be the graph obtained by
attaching an isomorph of $H$ to each of the $p$ nodes of $G$. Also, let
$\lambda_{1},\lambda_{2},\ldots,\lambda_{p}$ be the roots of
$B(G:x)$. Then
\begin{equation}
\label{15}
B(G(H);x)=\prod_{i=1}^{p}[B(H^{\triangle};x)-\lambda_{i}B(H_{-r};x)]\, .
\end{equation}
\end{lemma}

\noindent
Notice that $H^{\triangle}$ is misprinted in the original text (Corollary
2.2 of [19]) as $H$.

We can also derive a special corollary (see Corollary 2.3 in [19])
from Lemma 4, viz.:

\begin{lemma}\pnt
Let $0$ be a $k$-fold $(k\geq 1)$ root of $B(G;x)$. Then
$[B(H^{\triangle};x)]^{k}$ divides $B(G(H);x)$.
\end{lemma}

Here, we cannot help stating another lemma (see Lemma 6 in [20) that
generalizes Lemmas 4 and 5. First, denote by $H_{\lambda_{i}}$\, $(1\leq
i\leq p)$ the graph obtained by attaching a self-loop with the weight
$b_{r}=-\lambda_{i}$ to node $r$ of $H^{\triangle}$. It is not difficult
to establish that the expression in square brackets, in (15), is just
$B(H_{\lambda_{i}};x)$, which immediately affords us a derived result,
viz.:
\begin{lemma}\pnt
Let $G$ and $H$ be rooted graphs. Let $G(H)$ be the graph obtained by
attaching an isomorph of $H$ to each of the nodes of $G$. Also, let
$H_{\lambda_{i}}$\, $(1\leq i\leq p)$ be defined as above. Then
\begin{equation}
\label{16}
B(G(H);x)=\prod_{i=1}^{p}B(H_{\lambda_{i}};x)\, .
\end{equation}
\end{lemma}

One additional definition that will be employed below is this. 
Let $m_{1}(\lambda)$ and $m_{2}(\lambda)$ be the multiplicities of a
specific root $\lambda$ for polynomials $B(G_{1};x)$ and $B(G_{2};x)$,
respectively. We shall call the number $m(\lambda)=\min
(m_{1}(\lambda),m_{2}(\lambda))$ a {\em common multiplicity of an 
eigenvalue} $\lambda$ for the polynomials $B(G_{1};x)$ and
$B(G_{2};x)$.

Now we shall turn to deriving new results.

\section{Main results}

\hspace*{6mm}
The first our result will be complementary to Theorem 2:
\begin{theorem}\pnt
Let $G$ be a graph with $p$ vertices and
${\it\Gamma}=\{H_1,H_2,\ldots,H_p\}$ be a family of rooted graphs. Then
\begin{equation}
\label{17}
B(G({\it\Gamma});{\bf x};{\bf w})=\left[\prod_{i=1}^{p}B(L_i;{\bf x};{\bf
w})\right]\left[B(G^{\Box};{\bf x};{\bf
w})\left|_{x_i\rightarrow\frac{B(H_i;{\bf x};{\bf w})}{B(L_i;{\bf x};{\bf
w})}}\right.\right]\hspace{6mm}(1\leq i\leq p).
\end{equation}
\end{theorem}

\noindent
{\boldmath $Proof.$} Sketch it. Since the graphs $G$ and $H$ play a
symmetrical role in Lemma 1, one can rewrite (10) in an equivalent form as
follows
\begin{equation}
\label{18}
B(G\circ H;{\bf x};{\bf w})=B(G_{-r};{\bf x};{\bf w})\left[B(H;{\bf
x};{\bf w})\left|_{x_i\rightarrow\frac{B(G^{\triangle};{\bf x};{\bf
w})}{B(G_{-r};{\bf x};{\bf w})}}\right.\right]\hspace{6mm}(1\leq i\leq
p).
\end{equation}
Recall that Theorem 2 was proven in [19] by repetitively applying Lemma 
1 to $G({\it\Gamma})$. If we now use $p$ times (18) instead of (10),
we arrive at the result, wherein self-loops are (gradually)
removed from the core $G$ rather than from the the root $r_i$ of every
graph $H_i$ (in contrast to theorem 2). Therefore, denoting by
$G^{\Box}$ the core $G$ less all its self-loops (which is thus obtained),
we arrive at (17). Q.E.D. \hfill {\boldmath $\Box$}

Herein, we are interested in deriving a few corollaries of Theorem 7.
First, we shall state the following mate of Corollary 2.1:
\begin{corollary}\pnt
Let $T$ be a bipartite graph with the bipartition into $p_1$ and $p_2$
vertices, accordingly $(p_1\geq p_2;p_1+p_2=p)$. Let further
$T({\it\Gamma})$ be the restricted rooted product of graphs $T$ and
${\it\Gamma}$, wherein an isomorphic copy of a graph $H_1$ is attached to
every vertex of the first part of $T$ and an isomorph of another graph
$H_2$ is attached to every vertex of the second part of $T$ $(H_1,H_2\in
{\it\Gamma})$. Then
\begin{equation}
\label{19}
B(T({\it\Gamma});x)=\left[P(L_1;x)\right]^{p_1}\left[P(L_2;x)\right]^{p_2}
\left[P(T^{\Box};y_1,y_2)\left|_{y_i\rightarrow\frac{P(H_1;x)}{P(L_1;x)}}
\right.\right]\hspace{6mm}(i=1,2)\, ,
\end{equation}
where $T^{\Box}$ is the graph $T$ less all its self-loops; and $y_i$
$(i=1,2)$ simultaneously stands for all indexed $x$-variables belonging to
the vertices of the $i$-th part of $T$.
\end{corollary}

As an initial prerequisite to the next corollary, one can return to (13).
It gives an idea to expand $P(T^{\Box};y_1,y_2)$ of (19) as follows
\begin{equation}
\label{20}
P(T^{\Box};y_1,y_2)=\delta_{0}y_{1}^{p_1}y_{2}^{p_2}+\cdots+\delta_{k}
y_{1}^{p_1-k}y_{2}^{p_2-k}+\cdots+\delta_{p_2}y_{1}^{p_1-p_2}y_{2}^{0}
\hspace{6mm}(\delta_{0}=1;p_1\geq p_2;),
\end{equation}
where the adjunct powers of $y_1$ and $y_2$ decrease synchronously. The
matter is that all the $F$-covers (see (1)) that correspond to the circuit
polynomial $P(T^{\Box};y_1,y_2)$ of a bipartite loopless graph $T^{\Box}$
should consist just of cycles of even length (since only such cycles are
in it). Evidently, every $F$-cover always covers one and the same number
$k$ $(1\leq k\leq p)$ of vertices pertaining to both parts of $T^{\Box}$
(i.e., $k$ green vertices and $k$ red ones). Hence, it immediately
follows the above property of the powers in (20). The obtained expansion
affords the following corollary, of Theorem 7, accompanying Lemma 3:
\begin{lemma}\pnt
Let $P(T({\it\Gamma});x)$ be the circuit polynomial of the restricted
rooted product $T({\it\Gamma})$ as above. Then
\begin{equation}
\label{21}
P(T({\it\Gamma});x)=\sum_{k=0}^{p_2}\delta_{k}[P(H_1;x)]^{p_{1}-k}[P(H_2;x)]^
{p_{2}-k}[P(L_1;x)]^{k}[P(L_2;x)]^{k}\hspace{6mm}(p_1\geq p_2).
\end{equation}
\end{lemma}

As it was already discussed in Preliminaries, of our special interest are
instances $T(H)_1$ and $T(H)_2$ of $T({\it\Gamma})$ in which isomorphic
copies of an arbitrary graph $H$ are attached \underline{just} to the
vertices of either part of $T$ (while no attachment whatever is done to
the other part of it). This affords two complementary corollaries of Lemma
8.

\begin{corollary}\pnt
Let $T(H)_1$ be the restricted rooted product of a bipartite graph $T$ and
an arbitrary graph $H$ in which an isomorphic copy of $H$ is attached just
to every vertex of the first (greater) part of $T$. Then
\begin{equation}
\label{22}
P(T(H)_1;x)=\sum_{k=0}^{p_2}\delta_{k}(x+b_{2})^{p_{2}-k}[P(H;x)]^{p_{1}-k}
[P(H_{-r};x)]^{k}\hspace{6mm}(p_1\geq p_2),
\end{equation}
where $b_2$ is a common total weight of self-loops lying in each vertex of
the second (smaller) part of $T$.
\end{corollary}
\begin{corollary}\pnt
Let $T(H)_2$ be the restricted rooted product of a bipartite graph $T$ and
an arbitrary graph $H$ in which an isomorphic copy of $H$ is attached just
to every vertex of the second (smallerer) part of $T$. Then
\begin{equation}
\label{23}
P(T(H)_1;x)=\sum_{k=0}^{p_2}\delta_{k}(x+b_{1})^{p_{1}-k}[P(H;x)]^{p_{2}-k}
[P(H_{-r};x)]^{k}\hspace{6mm}(p_1\geq p_2),
\end{equation}
where $b_1$ is a common total weight of self-loops lying in each vertex of
the first (greater) part of $T$.
\end{corollary}

Now we recall that any circuit polynomial $P(T^{\Box};x)$ of a bipartite
graph $T^{\Box}$ without self-loops necessarily has at least
$(p_{1}-p_{2})$ zero eigenvalues (or roots), and together with every its
eigenvalue $\mu$ it also possesses an eigenvalue $-\mu$ (in particular see
Theorem 3.11 in [28]). In order to demonstrate this, one can substitute
$x$ for $y_1$ and $y_2$ on the R.H.S. of (20), which gives
\begin{equation}
\label{24}
P(T^{\Box};x)=\delta_{0}x^{p}+\delta_{1}x^{p-2}+\cdots+\delta_{p_{2}}x^
{p_{1}-p_{2}}\hspace{6mm}(p_1\geq p_2;p=p_{1}+p_{2}).
\end{equation}
It is immediately seen that $P(T^{\Box};x)$ is divisible by 
$x^{p_{1}-p_{2}}$ and, therefore, possesses at least $(p_{1}-p_{2})$ zero
eigenvalues. Further, all the powers of $x$ on the R.H.S. of (24) have one
and the same parity (either even or odd); and thereby the negative
$-\mu$ of every root $\mu$ of $P(T^{\Box};x)$ is also a root of it.
By this reason, we shall consider below only squares of the roots,
which, \underline{excluding} necessary $(p_{1}-p_{2})$\, 0's, will
comprise exactly $p_2$ $(p_2\leq p_1)$ not necessarily distinct numbers:
$\mu_{1}^{2},\mu_{2}^{2},\ldots,\mu_{p_2}^{2}$ (whose order does not
matter). This allows us to rewrite (24) as follows
\begin{equation}
\label{25}
P(T^{\Box};x)=x^{p_{1}-p_{2}}\prod_{i=1}^{p_2}(x^{2}-\mu_{i}^{2})\hspace{6mm}
(p_1\geq p_2).
\end{equation} 

Here, recall that, by definition, $\gamma_{k}$ is simultaneously the
coefficient of $x^{p-k}$ on the R.H.S. of (24) and of
$y_{1}^{p_{1}-k}y_{2}^{p_{2}-k}$ on the R.H.S. of (20); $1\leq k\leq
p_{2}\leq p_{1};p_{1}+p_{2}$. Among other things, this assures the reverse
passage from (24) to (20). But the R.H.S. of (24) is equal to the R.H.S.
of (25); therefore, we can legitimately rewrite (25) in two variables,
$y_1$ and $y_2$, as well:
\begin{equation}
\label{26}
P(T^{\Box};y_{1},y_{2})=y_{1}^{p_{1}-p_{2}}\prod_{i=1}^{p_{2}}(y_{1}y_{2}-
\mu_{i}^{2})\hspace{6mm}(p_{1}\geq p_{2}).
\end{equation}
In the present paper, (26) is an important requisite because it enables us
to state the following crucial sentence resembling Lemma 4:
\begin{lemma}\pnt
Let $P(T({\it\Gamma});x)$ be the circuit polynomial of the restricted
rooted product $T({\it\Gamma})$ of graphs $T$ and ${\it\Gamma}$ (see
above). Then
\begin{equation}
\label{27}
P(T({\it\Gamma});x)=\left[P(H_{1};x)\right]^{p_{1}-p_{1}}\prod_{i=1}^{p_2}\left
[P(H_{1};x)P(H_{2};x)-\mu_{i}^{2}P(L_{1};x)P(L_{2};x)\right]\hspace{6mm}(p_{1}
\geq p_{2}).
\end{equation}
\end{lemma}

\noindent
{\boldmath $Proof.$} Taking into account the definition of collective
variables $y_1$ and $y_2$ (instead of respective $x_i$'s in (17)) and
expressing $P(T^{\Box};{\bf x};{\bf w})$ in a specific form of the R.H.S
of (26), we can easily conclude that this statement is simply a corollary
of Theorem 7. Hence, the proof is immediate.\hfill {\boldmath $\Box$}

Note that interchanging the sorts of graphs $H_1$ and $H_2$ in
$T({\it\Gamma})$ (together with the weights $b_1$ and $b_2$) results in
another product $T({\it\Gamma})'$, which can be called, in the chemical
language, a {\em substitutional isomer} of $T({\it\Gamma})$. Under
$p_1=p_2$, the two substitutional isomers $T({\it\Gamma})$ and
$T({\it\Gamma})'$ distinguish only by the reciprocal fashion in which the
rooted graphs of sorts $H_1$ and $H_2$ are attached to the core $T$, in
them. Therefore, we shall call the last pair of substitutional isomers
{\em reciprocal rooted products}. This leads to the following corollary:
\begin{corollary}\pnt
Let $T$ be an equipartite bipartite graph $(p_1=p_2)$ and let
$T({\it\Gamma})$ and $T({\it\Gamma})'$ be the reciprocal (restricted) 
rooted products. Then
\begin{equation}
\label{28}
P(T({\it\Gamma});x)=P(T({\it\Gamma})';x).
\end{equation}
\end{corollary}

\noindent
{\boldmath $Proof.$} Indeed, under $p_1=p_2$, the indices 1 and 2 play
symmetrical roles on the R.H.S. of (27); consequently, the interchanging
of these indices cannot alter the result. Hence, the proof is
immediate.\hfill {\boldmath $\Box$}

Here, the chemical reader may recall that, in the reduced case, the role
of distinct chemical substituents, in a molecule, can be played by
heteroatoms. Or, in mathematical terms, a bipartite graph $T$ may possess
green vertices with the weight $b_1$ and red ones with the weight $b_2$;
interchanging $b_1$ and $b_2$ is just tantamount to the passage to the
reciprocal product (which is here simply a reweighted graph $T$, wherein
no (re-)attachments of any graphs $H_1$ and $H_2$ are made).

In a philosophical sense, it is very interesting that there exist
nonisomorphic graphs for which \underline{every specific} circuit
polynomial (characteristic, permanental, matching) should be equal. Here,
we recall that two graphs $G_1$ and $G_2$ are called {\em isospectral}
(or {\em cospectral}) (see [28]) if $P(G_1;x)=P(G_2;x)$; under this,
the type of the polynomial $P$ specifically depends on the context.
Studying isospectral graphs is an important aspect of the theory of graph
28]) and its application (see [11]). 

Regrettably, only in a descriptive form, the author dares to propose the 
simplest example of isospectral reciprocal products. First of all, note
that under $p_1=p_2$ the minimal bipartite graph in which the parts are
\underline{not equivalent} (on interchanging the colors of their vertices)
is a tree with 6 vertices. That tree is represented by a simple path
spanning five vertices, with the sixth vertex located at a free end of the
edge attached to the third (middle) vertex of it (see the tree 2.12 in
Table 2 of [28]). The chemical reader at once recognizes, in this tree, a
hydrogen-depleted graph of the carbon skeleton of 3-methyl-$n$-pentane
(where hydrogen atoms are not taken into account). We assume that the
green vertices are located on the longest 5-vertex chain of this tree at
sites 1, 3 and 5 while the red vertices are at sites 2, 4 and at the end
of the twig. Attaching an isomorph copy of an arbitrary graph $H_1$ (which
represents a chemical radical) to the green vertices and an isomorph of
another graph $H_2$ to the red vertices produces the product 
$T({\it\Gamma})$. The same process performed in a reciprocal fashion (when
sorts of graphs $H_1$ and $H_2$ are interchanged) does the reciprocal
rooted product $T({\it\Gamma})'$. Under this, the two restricted rooted
products are always isospectral; that is, the graphs $T({\it\Gamma})$ and 
$T({\it\Gamma})'$ should necessarily have one and the same circuit
polynomial. Further, let $H_1$ be a one-vertex graph $K_1$ (that is, no 
attachments should be made to the respective sites) and $H_2$ be a
two-vertex complete graph $K_2$ (i.e., an edge, or a twig). Then,
$T({\it\Gamma})$ (or the tree 2.74 in Table 2 [28]) is a
hydrogen-depleted graph of 2,4-dimethyl-3-ethyl-$n$-pentane while 
$T({\it\Gamma})'$ (or the tree 2.75 in Table 2 of [28]) is a
hydrogen-depleted graph of 4,$4'$-dimethyl-$n$-heptane.  As it follows
from Table 2 of the cited book, both of (molecular) graphs have the
same characteristic polynomial (of the adjacency matrix), viz.:
$$P(T({\it\Gamma});x)=P(T({\it\Gamma})';x)=x^{9}-8x^{7}+18x^{5}-12x^{3}\, 
,$$ with the roots (eigenvalues):                            
$\lambda_{1}=2.175;\lambda_{2}=1.414;\lambda_{3}=1.126;\lambda_{4}=
\lambda_{5}=\lambda_{6}=0;\lambda_{7}=-1.126;\lambda_{8}=-1.414;
\lambda_{9}=-2.175$. Possibly, [28] enables the reader to find other
instances of such graphs.

Now recall the definition of the restricted rooted products $T(H)_{1}$ and
$T(H)_{2}$, used in Corollaries 8.1 and 8.2, respectively. We can derive
similar corollaries of Lemma 9 as well:
\begin{corollary}\pnt
Let $P(T(H)_{1};x)$ be the circuit polynomial of the restricted rooted
product $T(H)_{1}$. Then
\begin{equation}
\label{29}
P(T(H)_1;x)=\left[P(H;x)\right]^{p_{1}-p_{2}}\prod_{i=1}^{p_{2}}
\left[(x+b_2)P(H;x)-\mu_{i}^{2}P(H_{-r};x)\right],
\end{equation}
where $b_2$ is the weight of every vertex of the second (smaller) part of
$T$.
\end{corollary}
\begin{corollary}\pnt
Let $P(T(H)_{2};x)$ be the circuit polynomial of the restricted rooted
product $T(H)_{2}$. Then
\begin{equation}
\label{29}
P(T(H)_1;x)=\left(x+b_1\right)^{p_{1}-p_{2}}\prod_{i=1}^{p_{2}}
\left[(x+b_1)P(H;x)-\mu_{i}^{2}P(H_{-r};x)\right],
\end{equation}
where $b_1$ is the weight of every vertex of the first (greater) part of
$T$.
\end{corollary}
Here, we remind that the weight $b_i$ $(i\leq i\leq p)$ of the $i$-th
vertex is a total weight of all self-loops lying in this vertex.

Lemma 9 and Corollaries 9.2 and 9.3 resemble Lemma 4. This list of
comparisons can be continued.
\begin{corollary}\pnt
Let 0 be an $s$-fold $(s\geq p_{1}-p_{2})$ root of the circuit polynomial
$P(T({\it\Gamma});x)$ of the restricted rooted product $T({\it\Gamma})$.
Then
$[P(H_1;x)]^{s}[p(H_2;x)]^{s-p_{1}+p_{2}}$ divides $P(T({\it\Gamma});x)$.
\end{corollary}

What roots $\mu$ are inherited by $(T({\it\Gamma});x)$ from the 
polynomials $P(H_1;x),P(H_2;x),P(L_1;x)$ and $P(L_2;x)$,  
involved in the above formulae, also depends on a common multiplicity
$m(\mu)$ of every root $\mu$ for four independent pairs
$(P(H_i;x);P(L_j;x))$ $(i,j\in \{1,2\})$ of these polynomials. Obviously,
if $m(\mu)\geq 1$ in at least one of the four cases the respective root
$\mu$ is also inherited by $P(T({\it\Gamma});x)$ (with the multiplicity
not less than $m(\mu)$). If $m(\mu)=0$ for all the four variants, it
(still) remains to employ Corollary 9.4, which works in a complementary
manner to the common-multiplicity principle. A more detailed investigation
is left here to the reader. However, we want to make yet some qualitative
remarks, addressed chiefly to the chemical audience. First, the degeneracy
of the roots of $P(T^{\Box};x)$ is beneficial for the multiplicity of the
roots of $P(T({\it\Gamma});x)$. Second, as well known, the root $\mu=0$ of
the characteristic polynomial of a molecular graph is unfavorable for the
stability of the molecule which is represented by it (see [11, 28]).
Nevertheless, the construction of the restricted rooted product can
harvest the same (or even more) benefits from the zero eigenvalue. That is
why the materials engineer must not a priori rule out the core graphs
$T^{\Box}$ with the eigenvalue(s) $\mu=0$. 

Now note that the factor $[P(H_1;x)P(H_2;x)-\mu_{i}^{2}P(L_1;x)P(L_2;x)]$
involved on the R.H.S. of (27) is, in its own right, the circuit polynomial 
$P(H_{\mu_{i}^{2}};x)$ of some derivative graph $H_{\mu_{i}^{2}}$. Here,
the graph $H_{\mu_{i}^{2}}$ is obtained by joining with the edge
$r_{1}r_{2}$ of the weight $\mu_{i}^{2}$ the roots $r_1$ and $r_2$ of
graphs $H_{1}$ and $H_2$, respectively; see the necessary general theory
in 20, 28]. Under this, it is worth recalling that the weight
$\mu_{i}^{2}$ of the edge $r_{1}r_{2}$ is, by definition, the product
$a'_{r_{1}r_{2}}a'_{r_{2}r_{1}}$ of the entries $a'_{r_{1}r_{2}}$ and
$a'_{r_{2}r_{1}}$ of the adjacency matrix $A'$ of the graph 
$H_{\mu_{i}^{2}}$. Our done preparation allows us to rewrite Lemma 9 into
the following equivalent form:

\begin{lemma}\pnt
Let $P(T({\it\Gamma});x)$ be the circuit polynomial of the restricted
rooted product $T({\it\Gamma})$ of graphs $T$ and ${\it\Gamma}$. Then
\begin{equation}
\label{31}
P(T({\it\Gamma});x)=[P(H_1;x)]^{p_{1}-p_{2}}\prod_{i=1}^{p_2}
P(H_{\mu_{i}^{2}};x)\hspace{6mm}(p_1\geq p_2),
\end{equation}
where $P(H_{\mu_{i}^{2}};x)$ is defined as above.
\end{lemma}

Lemma 10, derived for the circuit polynomial of the restricted rooted
product $T({\it\Gamma})$, resembles Lemma 6 for the $B$-polynomial of the
unrestricted rooted product $G(H)$ (see [19, 20]). Moreover, we notice
in passing that all the results above that involve the bipartite core
graph $T^{\Box}$ also hold good for the clique polynomial of
$T({\it\Gamma})$ (see [19]). 

Seeing the R.H.S. of (31), one can immediately conclude that there exists
such a similarity transformation of the adjacency matrix
$A(T({\it\Gamma}))$ of the restricted rooted product $T({\it\Gamma})$ that
block-diagonalizes it (to the form which is consistent with the R.H.S. of
(31)). The sense of the said can be better evaluated if one recalls how
difficult is, in general, to find whatever similarity transformation
simultaneously conserving all types of the circuit polynomial (of an
arbitrary matrix); usually, it is possible only for the characteristic 
polynomial while renders impossible for the permanental and/or matching
polynomial.

\section{Discussion}

\hspace*{6mm}
\noindent
\hspace{-7mm}1. By virtue of the results obtained herein, the circuit
polynomial $P(T({\it\Gamma});x)$ of the restricted rooted product
$T({\it\Gamma})$ can uniquely be reconstructed from the collection of
circuit polynomials $P(T^{\Box};x),P(H_1;x),P(H_2;x),P(L_1;x)$ and $P(L_2)$ 
of graphs $T^{\Box},H_1,H_2,L_1$ and $L_2$, consecutively. This
seems to be of use because, otherwise, it would be very difficult to
estimate the spectrum (of eigenvalues) of a complex target graph 
$T({\it\Gamma})$ (res. molecule), as is necessarily for creating
substances with given electronic and photonic properties. 

\noindent
2. An important role, in the present context, play degenerate roots of the
circuit polynomial $P(T({\it\Gamma});x)$. The degeneracy of eigenvalues,
first of all, tautologically means the possibility of filling some
(needed) energy levels, in a molecule or bulk material, to a rather
higher degree than it takes place for undegenerated eigenvalues. The said is 
profitable not only for electronic and photonic properties (of 
substances), as such, but can well be addressed to treating diverse
surfaces or substrates in order to passivate these or, on the contrary,
make more catalytically active, hydrophilic {\em etc.}. Under this, in
particular, a surface can acquire the properties of a chemical group
represented by a graph $H_1$ in cases when the latter multiply
contributes to the resulting spectrum of the system represented by the
graph $T({\it\Gamma})$ (see [20]).  

\noindent
3. A prospective trend, which is yet beyond present scope, is engineering
dendrimers (both graphs and molecules) that iteratively uses the
construction of the restricted rooted product $T({\it\Gamma})$. This can
essentially extend the capabilities of previous approaches (see [20]). 
Here, it is especially worth merging the methods of such work as [20]
with those obtained herein.

\section*{Acknowledgments}

\hspace*{6mm}
The author is indebted to the Referees for their attentive work with
his paper.

\bibliographystyle{plain}

\begin{thebibliography}{99}
   \bibitem{1} Issberner J., Moors R. and Vogtle F., Dendrimers: From
Generation and Functional Groups to Functions, {\em Angew. Chem. Int. Ed.
Engl.}, 1994, v. 33, no. 23/24, p. 2413--2420.    
   \bibitem{2} Archut A. and Vogtle F., Dendritic Molecules: Historic
Development and Future Applications, {\em Handb. Nanostruct. Mater. 
Nanotechnol}, 2000, v. 5, p. 333--374.
   \bibitem{3} Wang P.-W., Liu Y.-J., Devados C., Bharathi P. and Moore
J. S., Electroluminiscent Diods from a Single-Component Emitting
Layer of Dendritic Macromolecules, {\em Adv. Mater.}, 1996, v. 8, no. 3,
p. 237--241.  
   \bibitem{4} Bar-Haim A., Klafter J. and Kopelman R., Dendrimers
As Controlled Artificial Energy Antennae, {\em J. Am. Chem. Soc.}, 1997,
v. 119, p. 6197--6198.
   \bibitem{5} Bar-Haim A. and Klafter J., Dendrimers As Light-
Harvesting Antennae, {\em J. Lumin.}, 1998, v. 76, p. 197--200.
   \bibitem{6} Adronov A., Gilat S. L., Fr\'{e}chet J. M. J., Kaoru O.,
Neuwahl F. V. R., and Flemming G. R., Light Harvesting and Energy
Transfer in Laser-Dye-Labeled Poly(aryl ether)Dendrimers, {\em J. Am.
Chem. Soc.}, 2000, v. 122, p. 1175--1185. 
   \bibitem{7} Diudea M. V., Molecular Topology 21, Wiener Index of
Dendrimers, {\em Commun. Math. Comput. Chem. (MATCH)}, 1995, no. 32, p.
71--83.
   \bibitem{8} Diudea M. V. and Katona G., Molecular Topology of
Dendrimers, {\em Adv. Dendritic. Macromol.}, 1999, v. 4, p. 135--201. 
   \bibitem{9} Burioni R., Cassi D., Meccoli I. and Regina S., 
Tight-Binding Models on Branched Structures, {\em Phys. Rev. B: Condens.
Matter Mater.}, 2000, v. 61, no. 13, p. 8614--8617.              
   \bibitem{10} Yan J.-M. (Yen Ch.-M.), Symmetry Rules in the Graph
Theory of Molecular Orbitals, {\em Adv. Quantum Chem.} (Per-Olov
L\"{o}vdin, ed.), 1981, v. 13, p. 211--241.
   \bibitem{11} Dias J. R., {\em Molecular Orbital Calculations Using
Chemical Graph Theory}, Springer-Verlag, Berlin, 1993.
   \bibitem{12} Dias J. R., Techniques in Facile Calculation of
Molecular Orbital Parameters and Related Conceptualizations---Molecular
Orbital Functional Groups, {\em J. Mol. Struct. (Theochem)}, 1997, no.
417, p. 49--67.
   \bibitem{13} Dias J. R., From Small Molecules to Infinitely Large
$\pi$-Electron Networks---Strongly Subspectral Molecular Systems, {\em Z.
Naturforsch}, 1998, v. 53a, p. 909--918.
   \bibitem{14} Dias J. R., Analysis of $\pi$-Electronic Structures of
Small Alternant Hydrocarbons to Infinitely Large Polymeric Strips. The
Aufbau Principle and End-Group Effects, {\em Int. J. Quantum Chem.}, 1999,
v. 74, p. 721--724.
   \bibitem{15} Dias J. R., Two-Dimensional Arrays in the Analysis of
Trends in Series of Molecules: Strongly Subspectral Molecular Graphs,
Formula Periodic Tables, and Number of Resonance Structures, {\em J. Chem.
Inf. Comput. Sci}, 2000, v. 40, no. 3, p. 810--815.
   \bibitem{16} Dias J. R., Strongly Subspectral Series Containing
Cyclobutadiene Moiety, {\em Croat. Chem. Acta}, 2000, v. 73, no. 2, p.
405--415.
   \bibitem{17} Rosenfeld V. R., Endomorphisms of a Weighted Molecular
Graph and Its Spectrum, {\em Commun. Math. Comput. Chem. (MATCH)}, 1999,
no. 40, p. 203--214.
   \bibitem{18} Rosenfeld V. R., Some Spectral Properties of the
Arc-Graph, {\em Commun. Math. Comput. Chem. (MATCH)}, 2001, no. 43, p.
41--48. 
   \bibitem{19} Farrell E. J. and Rosenfeld V. R., Block and
Articulation
Node Polynomials of the Generalized Rooted Product of Graphs, {\em Jour.
of Mathematical Sciences} (India), 2000, v. 11, no. 1, p. 35--47.
   \bibitem{20} Rosenfeld V. R. and Diudea M. V., The Block Polynomials
and Block Spectra of Dendrimers, {\em Internet Electron. J. Mol. Des.},
2002, v. 1, no. 3, p. 142--156.
   \bibitem{21} Farrell E. J., On a General Class of Graph Polynomials,
{\em J. Comb. Theory B}, 1979, v. 26, no. 1, p. 111--122.
   \bibitem{22} Rosenfeld V. R. and Gutman I. M., A Novel Approach to
Graph Polynomials, {\em Commun. Math. Comput. Chem (MATCH)}, 1989, no. 24,
p. 191--199. 
   \bibitem{23} Rosenfeld V. R. and Gutman I. M., On the Graph
Polynomials
of a Weighted Graph, {\em Coll. Sci. Papers. Fac. Kragujevac}, 1991, v.
12, p. 49--57.
   \bibitem{24} Farrell E. J., On a Class of Polynomials Obtained from
Circuits in a Graph and Its Application to Characteristic Polynomials of
Graphs, {\em Discrete Math.}, 1979, v. 25, p. 121--133. 
   \bibitem{25} Farrell E. J. and Grell J. C., The Circuit Polynomial
and Its Relation to Other Polynomials, {\em Carrib. J. Math.}, 1982, v.
2, no. 1/2, p. 15--24. 
   \bibitem{26} Kerber A., {\em Algebraic Combinatorics via Finite
Group Actions}, Wissenschaftsverlag, Manheim, Wein, Z\"{u}rich, 1991.
   \bibitem{27} Kerber A., {\em Applied Finite Group Actions}, Springer
Verlag, Berlin, Heidelberg, New York, London, Raris, Tokyo, Hong Kong,
Barcelona, Budapest, 1999.
   \bibitem{28} Cvetkovi\'{c} D. M., Doob M. and Sachs H., {\em Spectra
of Graphs: Theory and Application}, Academic Press, New York, 1980.
   \bibitem{29} Cvetkovi\'{c} D. M., Doob M., Gutman I. M. and
Torga\v{s}ev A., {\em Recent Results in the Theory of Graph Spectra},
North-Holland, Amsterdam, 1988.
   \bibitem{30} Godsil C. and McKay B. A., A New Graph Product and Its
Spectrum, {\em Bull. Austral. Math. Soc.}, 1978, v. 18. p. 21--28.
   \bibitem{31} Farrell E. J., An Introduction to $F$-Graphs, a
Graph-Theoretic Representation of Natural Numbers, {\em Internat. J. Math.
and Math. Sci.}, 1992, v. 15, no. 2, p. 313--318.
   \bibitem{32} Rosenfeld V. R. and Rosenfeld Victor R., Groupoids and
Classification of Polymerization Reactions, in {\em The Use of Computers
in Spectroscopy and Chemical Research, Novosibirsk, 6--8th September 1983.
Theses of the All-Union Conference}, Novosibirsk, 1983, p. 195--196
(in Russian).
   \bibitem{33} Rosenfeld V. R., An Algebraic Model of Closed Loops in
Proteins, {\em Commun. Math. Comput. Chem. (MATCH)}, submitted.
\end{thebibliography}

\end{document}